\documentclass[a4paper,fleqn,11pt]{amsart}

\usepackage{amssymb, comment}

\usepackage{mathrsfs}

\usepackage{graphicx}
\graphicspath{{./},{./images/}}

\usepackage{graphicx}
\usepackage{multirow}
\usepackage[colorlinks=true]{hyperref}
\numberwithin{equation}{section}

\newcommand{\vb}{\vspace{3mm}}

\def \approxd {\,{\buildrel d \over \approxd}\,}

\theoremstyle{remark}

\renewcommand{\proofname}{\noindent {\bf Proof}}

\begin{document}
\title[Congestion analysis of unsignalized intersections]{Congestion analysis of unsignalized intersections:\\The impact of impatience and Markov platooning}
\author{Abhishek$^1$}\thanks{$1$ Korteweg-de Vries Institute for Mathematics, University of Amsterdam, Amsterdam, The Netherlands
({\tt \{abhishek,m.r.h.mandjes,nunezqueija\}@uva.nl})}

\author{Marko Boon$^2$}\thanks{
$^2$Department of Mathematics and Computer Science, Eindhoven
University of Technology, P.O. Box 513, 5600 MB Eindhoven, The
Netherlands ({\tt m.a.a.boon@tue.nl})}


\author{Michel Mandjes$^1$}\thanks{
}

\author{Rudesindo N\'u\~nez-Queija $^1$}\thanks{
}

\date{\today}


\begin{abstract}
This paper considers an unsignalized intersection used by two traffic streams. A stream of cars is using a primary road, and has priority over the other
stream. Cars belonging to the latter stream cross the primary road if the gaps between two subsequent cars on the primary road are larger than their critical headways. A question that naturally arises relates to the capacity of the secondary road: given the arrival pattern of cars on the primary road, what is the maximum arrival rate of low-priority cars that can be sustained? This paper addresses this issue by considering a compact model that sheds light on the dynamics of the considered unsignalized intersection. The model, which is of a queueing-theoretic nature, reveals interesting insights into the impact of the user behavior on stability.

\noindent The contributions of this paper are threefold. First, we obtain new results for the aforementioned model that includes driver impatience. Secondly, we reveal some surprising aspects that have remained unobserved in the existing literature so far, many of which are caused by the fact that the capacity of the minor road cannot be expressed in terms of the \emph{mean} gap size; instead more detailed characteristics of the critical headway distribution play a crucial role. The third contribution is the introduction of a new form of bunching on the main road, called {\it Markov platooning}. The tractability of this model allows us to study the impact of various platoon formations on the main road on the capacity of the minor road.

\end{abstract}

\maketitle


\section{Introduction}
A common element in road traffic networks is that of an unsignalized intersection that is used by two traffic streams which have different priorities. In the first place there is a high-priority class that consists of cars that use a major (or primary) road. These cars pass the intersection according to some inherently random process. Having priority, they do so without observing the low-priority stream. Cars of the low-priority stream, which use a minor (or secondary) road, however, only cross when the duration (in time) of a gap between two subsequent cars passing by on the main road is sufficiently large, i.e., larger than a (possibly car-specific) threshold $T$.

As the high-priority cars on the primary road do not experience any interference from the low-priority cars, the system's performance is fully determined by the characteristics of the queue of low-priority cars on the secondary road. A first topic of interest is the capacity of this secondary road, which is defined as the maximum possible number of departures (per time unit) of vehicles on this road. Heidemann and Wegmann \cite{heidemann97} show that this definition implies that the capacity can be expressed in terms of the stability of the corresponding queue: for what arrival rate of low-priority cars can it be guaranteed that the queue remains bounded? The answer to this question evidently depends on the distribution of the gaps between subsequent cars on the primary road. In particular, the capacity of the minor road is greatly influenced by the clustering of vehicles in platoons on the main road. In addition, specific features of the low-priority car drivers play a crucial role, in terms of the way that individual car drivers choose their critical headways. In the existing literature, various models have been studied, the simplest variant being that all low-priority drivers use the same deterministic critical headway $T$ \cite{guo}. A second, more realistic, model allows different values of $T$ according to some probabilistic distribution \cite{guo2014,wu2012}, where $T$ is resampled for any new attempt at crossing the main road. The randomness captures the heterogeneity in the preferences (and driving styles) of the low-priority car drivers. A further refinement is a model in which different drivers have different thresholds $T$, but in which each driver persistently uses a single driver-specific value of $T$ for all attempts.
In this paper we will investigate these issues, building on and extending our preliminary results in the short paper~\cite{abhishekcomsnets2016}.
For consistency with that work, we will denote the three behavior types described above by B$_1$, B$_2$, and B$_3$.

\vspace{2mm}

Various aspects of gap acceptance models have been studied before. The main applications concern  unsignalized intersections (e.g.\ \cite{catchpoleplank,cheng,heidemann97,tanner62}), pedestrian crossings (e.g.\ \cite{mayne,tanner51,wei}), and freeways (e.g.\ \cite{drew2,drew1}). Although the gap acceptance process in these three application areas exhibits similar features, the queueing aspects are fundamentally different. In this paper, we focus on motorized vehicles, but all results regarding the capacity of the minor road can be applied to pedestrian crossings or freeway merging. Heidemann and Wegmann \cite{heidemann97} give an excellent overview of the existing results in gap acceptance theory, including the three types of user behavior that were discussed above.

\vspace{2mm}

Several relevant aspects have not yet been incorporated in previously investigated models.
The main objective of this paper is to further enhance this class of models by extending the framework in~\cite{abhishekcomsnets2016} with platooning on the main road.
Our work contains the following three contributions.
\begin{itemize}
\item[$\circ$]
First, we show how to incorporate  \emph{impatience} of the drivers that are waiting to cross the major road. This phenomenon, which is indeed encountered in practice \cite{abouhenaidy94}, has been studied before in e.g.\ \cite{drew2, drew1, weissmaradudin}, but (to the best of our knowledge) not yet in the context of models B$_2$ and B$_3$, where randomness is encountered in the critical headway $T$.
\item[$\circ$]The second contribution concerns a number of surprising aspects that have remained unobserved in the existing literature so far. We show that  the capacities that correspond with the three different types of driver behavior that we introduced above, are strictly ordered: B$_2$ has the largest capacity, then B$_1$, and the capacity of B$_3$ is the smallest (with the mean critical headway of models B$_2$ and B$_3$ chosen equal to the deterministic critical headway of model B$_1$). As it turns out, the capacity can {\it not} be given in terms of the mean quantity ${\mathbb E}[T]$, but more precise distributional information of the random variable $T$ is needed. Perhaps counterintuitively, when comparing two gap time distributions $T_1$ and $T_2$ one could for instance encounter situations in which ${\mathbb E}[T_1]< {\mathbb E}[T_2]$, but in which still the capacity under $T_1$ is smaller than the one under $T_2$.
\item[$\circ$]The third contribution is the introduction of a new model for vehicle clustering on the main road, which we will refer to as {\it  Markov platooning} throughout this paper. The tractability of this model allows us to study the impact of various platoon formations \cite{gaur2001,jia2016,li2017} on the main road on the capacity of the minor road. Platoon forming has also been studied in the existing literature on gap acceptance models before. The most common models that include clustering on the major road are so-called gap-block models. In these models, vehicles tend to form platoons, most commonly arriving according to Poisson processes.  The lengths of these platoons are i.i.d. random variables with general distributions, which can be chosen carefully to mimic real-life clustering behavior. Tanner \cite{tanner62} considers a model where platoon lengths are distributed as the busy period of a single-server queue. Wegmann \cite{wegmann1991} and Wu \cite{wu2001} analyze the capacity under even less restrictive assumptions. However, all of these models assume no (or a very weak form of) dependence between successive block sizes and gap sizes. By introducing Markov platooning, an arrival process based on Markov modulation, we allow for a more refined way of bunching on the major road that includes dependence between successive gap sizes.
\end{itemize}

This paper is structured as follows. In the next section, we describe in more detail the variations of the gap acceptance model, including the aforementioned types of gap acceptance behavior, impatience, and platooning on the major road. In Section \ref{sect:analysis} we analyze queue lengths and delays in the standard model and present numerical results for several practical examples, focusing on some surprising, paradoxical features that one might encounter. In Sections \ref{sect:impatience} and \ref{sect:markovplatooning}, we study the impact of impatience and Markov platooning on traffic congestion on the minor road. In these sections we also present numerical results, exhibiting interesting features of the model variations. Section \ref{sect:conclusions} concludes the paper.

\section{Preliminaries}\label{sect:model}

\subsection{Arrival process.}
\begin{figure}[!ht]
\begin{center}
\includegraphics[width=0.7\linewidth]{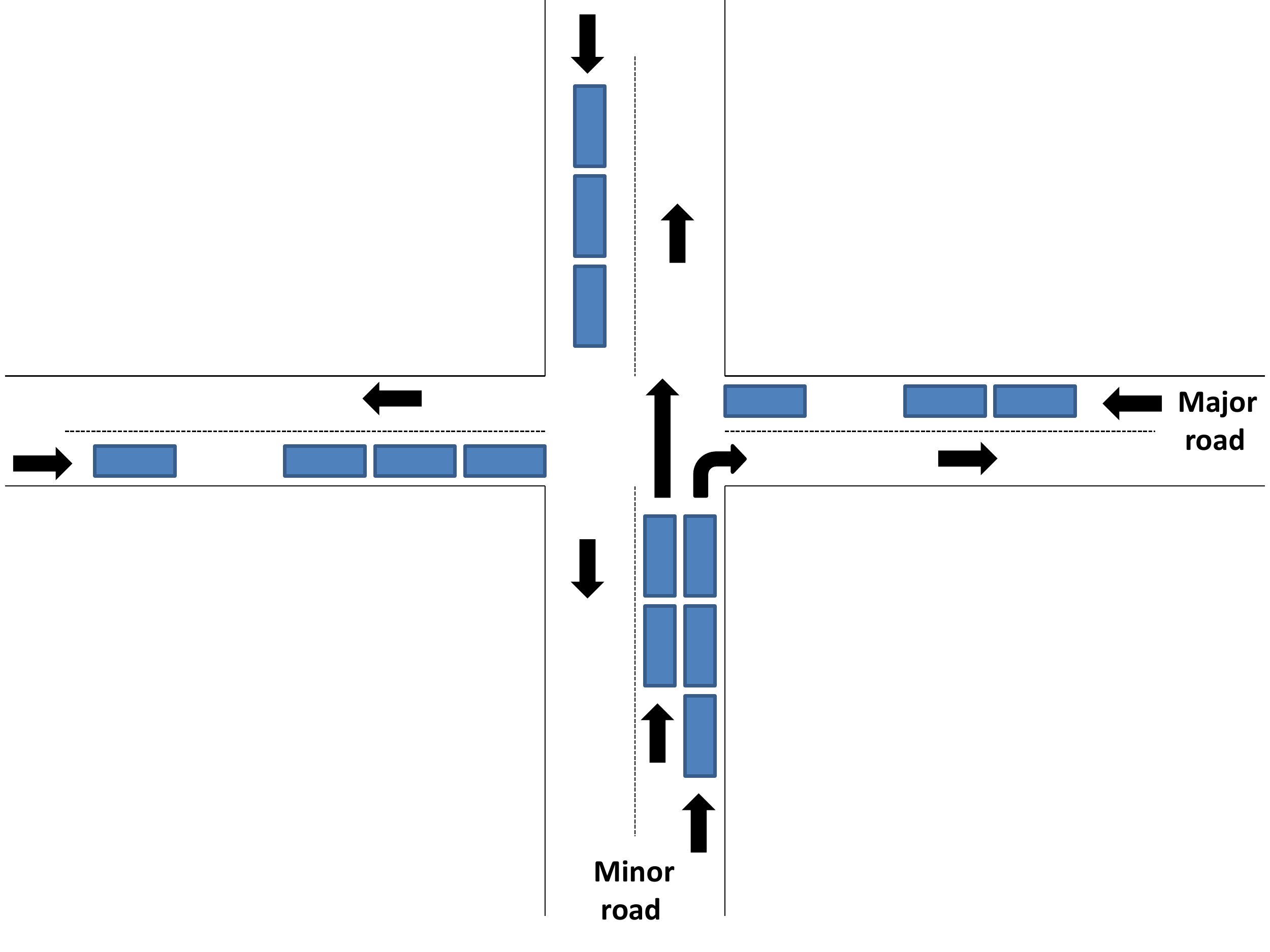}
\end{center}
\caption{An example of a situation that can be analyzed using the model in this paper.}
\label{fig:intersection}
\end{figure}

The situation analyzed in this paper is depicted in Figure~\ref{fig:intersection}. We consider an intersection used by two traffic streams, both of which wishing to cross the intersection. There are two priorities: the cars on the major road have priority over cars on the minor road (and hence do not notice the presence of the minor road). The low-priority cars on the minor road cross the intersection as soon as the gap between two subsequent high-priority cars has a duration larger than $T$, commonly referred to as the \emph{critical headway}.

Cars on the minor road arrive according to a Poisson process with rate $\lambda$. In this paper we distinguish between two types of arrival processes on the major road. The first arrival process for the high-priority car drivers, which we consider in Sections \ref{sect:analysis} and \ref{sect:impatience}, is a classical Poisson process with intensity $q$, meaning that the inter-arrival times between any pair of subsequent cars are exponentially distributed with mean $1/q$. The second arrival process is a generalization of the Poisson process, viz.\ the Markov modulated Poisson process (MMPP). The MMPP, which will be discussed in greater detail in Section \ref{sect:markovplatooning}, is a well-studied arrival process which is generally used to model dependencies between inter-arrival times. In an MMPP, at time $t$ the time till the next arrival is exponentially distributed with mean $1/q_i$ if an independently evolving Markov process (usually referred to as the {\it background process}) is in state $i$ at time $t$.  The flexibility of the MMPP allows us to vary the inter-arrival times in such a way, that we can create platoons, single arrivals, or combinations thereof.

\subsection{Gap acceptance behavior.}
We have not yet exactly defined the criterion by which the low-priority cars decide to cross. In this paper we distinguish three types of `behavior' when making this decision.
\begin{itemize}
\item[B$_1$]
The first model is the most simplistic: the critical headway $T$ is deterministic, and uniform across all low-priority car drivers.
\item[B$_2$]  Clearly B$_1$ lacks realism, in that there will be a substantial level of heterogeneity in terms of driving behavior: one could expect a broad range of `preferences', ranging from very defensive to very reckless drivers. In B$_2$ this is modeled by the car driver at the front end of the queue {\it resampling} $T$ (from a given distribution) at any new attempt (where an `attempt' amounts to comparing this sampled $T$ to the gap between the two subsequent cars that he is currently observing).
\item[B$_3$] In the third model an alternative type of driver behavior is assumed. More specifically,  it reflects {\em persistent} differences between drivers, in that each driver selects a random value of $T$, but then sticks to that same value for all attempts, rather than resampling these.
\end{itemize}

\subsection{Impatience.}
For each of the aforementioned behavior types, we also consider a variant that includes impatience. With impatience, the critical headway decreases after each failed attempt, reflecting the impatience of drivers, resulting in the willingness to accept smaller and smaller gaps. In more detail, we define a critical headway $T_j$ for the $j$-th attempt to enter the main road ($j=1,2,\dots$). Note that, depending on the distributions of $T_1, T_2, \dots$, in model B$_2$ situations might occur where $T_{i+1}>T_i$, despite $T_{i+1}$ being stochastically smaller than $T_i$. This is a typical feature of the model with resampling.
Exact details regarding the manner in which impatience is incorporated will be given in the next section.

A few remarks are in place here. In the first place, above we positioned this setup in the context of an unsignalized intersection, but various other applications could be envisioned. One of these could correspond to the situation in which the low-priority cars have to merge with the stream of high-priority cars (e.g. from a ramp or a roundabout). Also in the context of pedestrians crossing a road, the model can be used. We also stress that in the case the primary road actually consists of two lanes that have to be crossed (without a central reservation), with cars arriving (potentially in opposite directions) at Poisson rates (say) $q^{\leftarrow}$ and $q^{\rightarrow}$, our model applies as well, as an immediate consequence of the fact that  the superposition of two Poisson processes is once again a Poisson process with the parameter $q:=q^{\leftarrow}+q^{\rightarrow}$; see also the discussion in \cite[Section 5]{wu2001}.

\section{New insights for the classical model}\label{sect:analysis}

In this paper we analyze the three models relying on queueing-theoretic techniques.
In this section we consider the classical setting with Poisson arrivals and no impatience.
Since many results for the variants without impatience have been known in the existing literature (see, for example, Heidemann and Wegmann \cite{heidemann97} for an overview), we will mainly focus on the additional insights that can be obtained by carefully studying the formulas for the capacity of the minor road under different circumstances, which turns out to lead to a few interesting new insights.

We start our exposition by introducing some notation. In the first place, we let $X_n$ denote the number of cars in the queue on the minor road when ({\it right after}, that is) the $n$-th low-priority car crosses the primary road; in addition, $T_n^\#$ is the time that this happens.
We let $Y_n:=T_n^\#-T_{n-1}^\#$ denote the inter departure time between the $(n-1)$-st and $n$-th car from the secondary road. 
It is well-known that the process $\{X_n, n=1,2,...\}$ has the dynamics of a standard single-server queue with Poisson arrivals and general service times, in Kendall's famous notation also known as the M/G/1 queue (although in this context it is more common to refer to the model as an M/G2/1 queue). The dynamics of the merging process is fully captured in the distribution of the service times, $Y_1, Y_2, \dots$, which will have their specific form for each of the models B$_1$ up to B$_3$. As a consequence, we have that $X_n$ has a stationary distribution which is uniquely characterized through its  probability generating function (directly following from  the celebrated Pollaczek--Khinchine formula).
For a formal derivation of some of the expressions in this section, we refer to our earlier paper \cite{abhishekcomsnets2016}. In the current paper we do not focus on queue lengths or delays, but we only focus on the impact of the three types of the driver's behavior on the `capacity' of the secondary road; here `capacity' is defined as the maximum arrival rate $\lambda$ such that the corresponding queue does not explode. A standard result from queueing theory is that for the M/G/1 queue the stability condition is $\rho:=\lambda\,{\mathbb E}[Y]<1$. As a consequence, the capacity of the minor road, denoted by $\bar{\lambda}$, can be determined for each of the models, with or without impatience (see Heidemann and Wegmann \cite{heidemann97}):
\begin{equation}\label{fourc}
\bar{\lambda}=\frac{1}{{\mathbb E}{[Y]}},
\end{equation}
where $\mathbb {E}{[Y]}$ depends on the driver behavior, as explained before. In this section we denote the capacity for model B$_i$ by $\bar{\lambda}_i$, for $i=1,2,3$. Although the expressions below can also be found in, for example, Heidemann and Wegmann \cite{heidemann97}, we add some new observations regarding the capacities. Since the results below are special cases of the variants with impatience, they can also be obtained using the results in Section \ref{sect:impatience}.

\vspace{2mm}
\noindent
\textbf{\boldmath B$_1$ {(constant gap):}} Every driver on the minor road needs the same constant critical headway $T$ for evert attempt to enter the main road.
\begin{align*}
{\mathbb E} [e^{-sY}] = \frac{(s+q)e^{-(s+q)T}}{s+qe^{-(s+q)T}},\quad
{\mathbb E}[Y] = \frac{e^{qT}-1}{q},\quad
\bar{\lambda}_1:=\frac{q}{e^{qT}-1}.
\end{align*}

\vspace{2mm}

\noindent
\textbf{\boldmath B$_2$ (sampling   per attempt):} With this behavior type, which is also sometimes referred to as ``inconsistent behavior'', every car driver samples a random $T$ for each new `attempt' (where `attempt' corresponds to comparing the resulting $T$ with the gap between two subsequent cars on the major road).
\begin{align*}
{\mathbb E} [e^{-sY}] = \frac{(s+q){\mathbb E}[e^{-(s+q)T}]}{s+q\,{\mathbb E}[e^{-(s+q)T}]},\quad
\mathbb{E}[Y]=\frac{1-\mathbb{E}[e^{-qT}]}{q\mathbb{E}[e^{-qT}]}, \quad\bar{\lambda}_2=\frac{q}{({\mathbb E}[e^{-qT}])^{-1}-1}. \label{six}
\end{align*}

\noindent
\textbf{\boldmath B$_3$ (sampling per driver):}   In this variant, sometimes referred to as ``consistent behavior'', every car driver samples a random $T$ at his first attempt. This (random) value will be used consistently for each new attempt by this driver.

\begin{align*}
{\mathbb E} [e^{-sY}] = {\mathbb{E}}\left[\frac{(s+q)e^{-(s+q)T}}{s+qe^{-(s+q)T}}\right],\quad
{\mathbb E}[Y] = \frac{\mathbb{E}[e^{qT}]-1}{q}, \quad \quad\bar{\lambda}_3=\frac{q}{\mathbb{E}[e^{qT}]-1}.
\end{align*}

Importantly, it is here tacitly assumed that the moment generating function $\mathbb{E}[e^{qT}]$ of $T$ exists. A consequence that has not received much attention in the existing literature, is that it also means that in case $T$
has a polynomially decaying tail distribution (i.e., ${\mathbb P}(T>t) \approx C\,t^{-\beta}$ for some $C,\beta>0$ and $t$ large) {\it the queue at the secondary road is never stable}. The reason is that for this type of distributions it is relatively likely that an extremely large $T$ is drawn, such that it takes very long before the car can cross the intersection (such that in the mean time the low-priority queue has built up significantly).

In fact, also for certain light-tailed distributions we find that B$_3$ has an undesirable impact on the capacity. Take, for example, $T$ exponentially distributed with parameter $\alpha$. In this case, we have
\[
\mathbb{E}[Y]=\begin{cases}
1/(\alpha-q), & \quad q < \alpha,\\
\infty & \quad q \geq \alpha,
\end{cases}
\]
implying that the capacity of the minor street drops to zero when $q \geq 1/\mathbb{E}[T]$. Actually, the situation might be even worse than it seems, because it can be shown that  $\mathbb{E}[Y^k]=\infty$ if $q \geq \alpha/k$, for $k=1,2,\dots$. As a consequence, when $\alpha > q \geq \alpha/2$ the capacity is positive, but the mean queue length and the mean delay at the minor road grow beyond any bound.\\

Another interesting observation, is that the arrival rates $\bar{\lambda}_1, \bar{\lambda}_2$ and $\bar{\lambda}_3$ obey the ordering
\[\bar\lambda_2\geqslant \bar\lambda_1\geqslant \bar\lambda_3\]
(where in B$_1$ we have chosen $T$ equal to the mean ${\mathbb E}[T]$ used in the other variants). This is an immediate consequence of Jensen's inequality, as we show now. To compare $\bar{\lambda}_1$ and $ \bar{\lambda}_3$ realize that Jensen's inequality implies
\[ \frac{1}{q} (\mathbb{E}[e^{qT}]-1)\geqslant \frac{1}{q}(e^{q\mathbb{E}[T]}-1),\]
which directly entails $\bar{\lambda}_3 \leq \bar{\lambda}_1$. Along the same lines, again appealing to Jensen's inequality,
$
\mathbb{E}[e^{-qT}]\geqslant e^{-q\mathbb{E}[T]},$ and hence
$\bar{\lambda}_2 \geq \bar{\lambda}_1$.

\vspace{2mm}

We conclude this section by stating a number of general observations and illustrative numerical examples. In the first place, the above closed-form expressions show that {the stability conditions (and hence the capacities) depend on the full distribution of $T$}, as opposed to just the mean value ${\mathbb E}[T]$.

In many stochastic dynamic systems introducing variability leads to a degradation of  the performance. The fact that $\bar\lambda_2 \geqslant \bar\lambda_1$ indicates that in this case this `folk theorem' does not apply: the fact that one resamples $T$ often (every driver selects a new value for each new attempt) actually {\it increases} the capacity of the low-priority road.

\subsection{Example 1: ordering of the capacities}

\begin{figure}[!ht]
\parbox{0.49\textwidth}{\centering
\includegraphics[width=\linewidth]{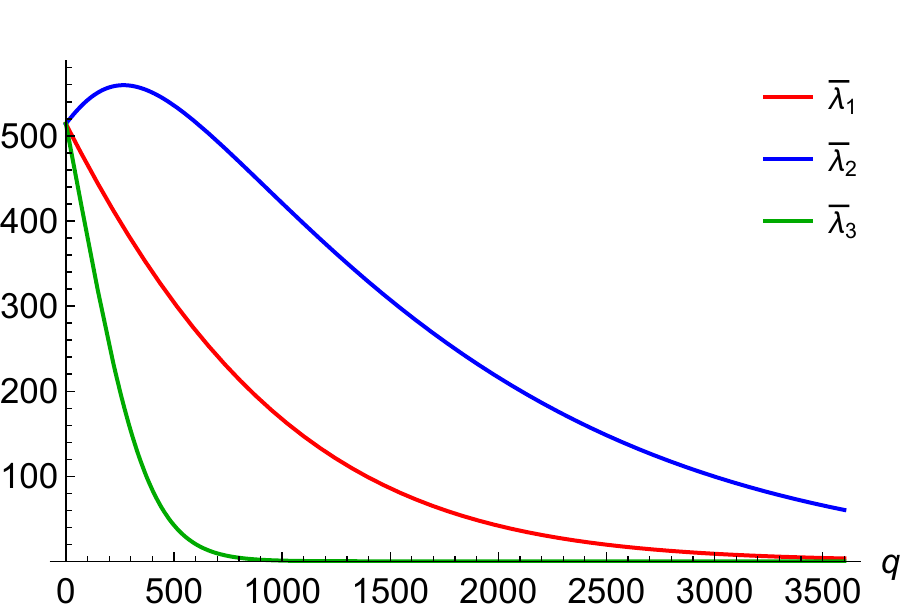}\\
(a) Capacity  vs. $q$ (veh/h)
} \hfill
\parbox{0.49\textwidth}{\centering
\includegraphics[width=\linewidth]{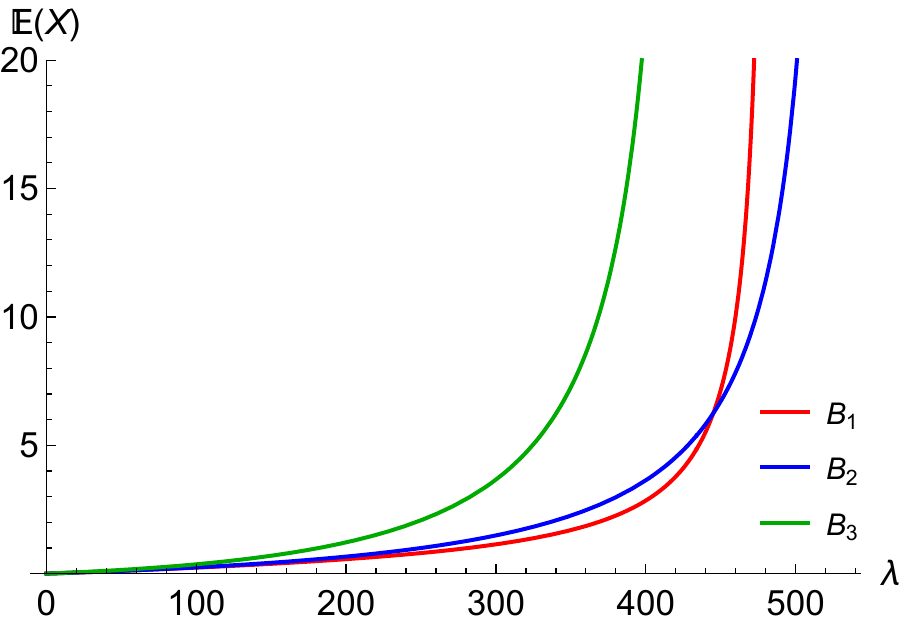}\\
(b) Mean queue lengths  vs. $\lambda$ (veh/h)
} \hfill
\caption{The capacities (veh/h) and mean queue lengths in Example 1.}
\label{fig:example1capacity}
\end{figure}

In this example we illustrate the impact of driver behavior on the capacity of the system and on the queue lengths. In particular, we compare the following three scenarios (corresponding to the three behavior types), with the parameters chosen such that the system exhibits interesting features:
\begin{itemize}
\item[(1)] All drivers search for a gap between consecutive cars on the major road, that is at least $7$ seconds long.
\item[(2)] A driver on the minor street, waiting for a suitable gap on the major street, will sample a new (random) critical headway every time a car passes on the major street. With probability $9/10$ this critical headway is $4$ seconds, and with probability $1/10$ it is exactly $34$ seconds. Note that the expected critical headway is $0.9\times4+0.1\times34=7$ seconds, ensuring a fair comparison between this scenario and the previous scenario.
\item[(3)] In this scenario we distinguish between slow and fast traffic. We assume that 90\% of all drivers on the minor road need a gap of (at least) $4$ seconds. The other 10\% need at least $34$ seconds.
\end{itemize}
Figure\ \ref{fig:example1capacity}(a) depicts the capacity (veh/h) of the minor street as a function of $q$, the flow rate on the main road (veh/h). The relation $\bar\lambda_2\geqslant\bar\lambda_1\geqslant\bar\lambda_3$ is clearly visible. An interesting aspect, exhibited by Figure \ref{fig:example1capacity}(a), is that the capacity for Model~B$_2$ is not always decreasing in $q$, but we return to this topic in the next numerical example.
In this example we focus on the comparison of the mean queue lengths resulting from the different types of driver's behavior. Although the capacities for B$_1$, B$_2$, and B$_3$ are strictly ordered, $\lambda_2 \geqslant \lambda_1 \geqslant \lambda_3$, this is not necessarily true for the mean queue lengths on the minor street. To illustrate this, we fix the traffic flow on the major road at $q=60$ vehicles per hour. The mean queue length on the minor road, as a function of $\lambda$, is depicted in Figure\ \ref{fig:example1capacity}(b). The figure displays a paradoxical situation, where the mean queue length corresponding to B$_2$ (resampling) is \emph{higher} than with B$_1$ (constant) for $71.2 < \lambda < 445.1$. Indeed, for $\lambda < 71.2$, which is hardly visible but definitely computable, B$_2$ also has a smaller mean queue length than B$_1$. It can be shown that, for the given distribution of $T$, this paradox only takes place when $q < 124.6$ vehicles per hour. Note that, due to Little's law, the mean delays exhibit the same behavior.

\subsection{Example 2: the impact of resampling}

\begin{figure}[!ht]
\includegraphics[width=0.6\linewidth]{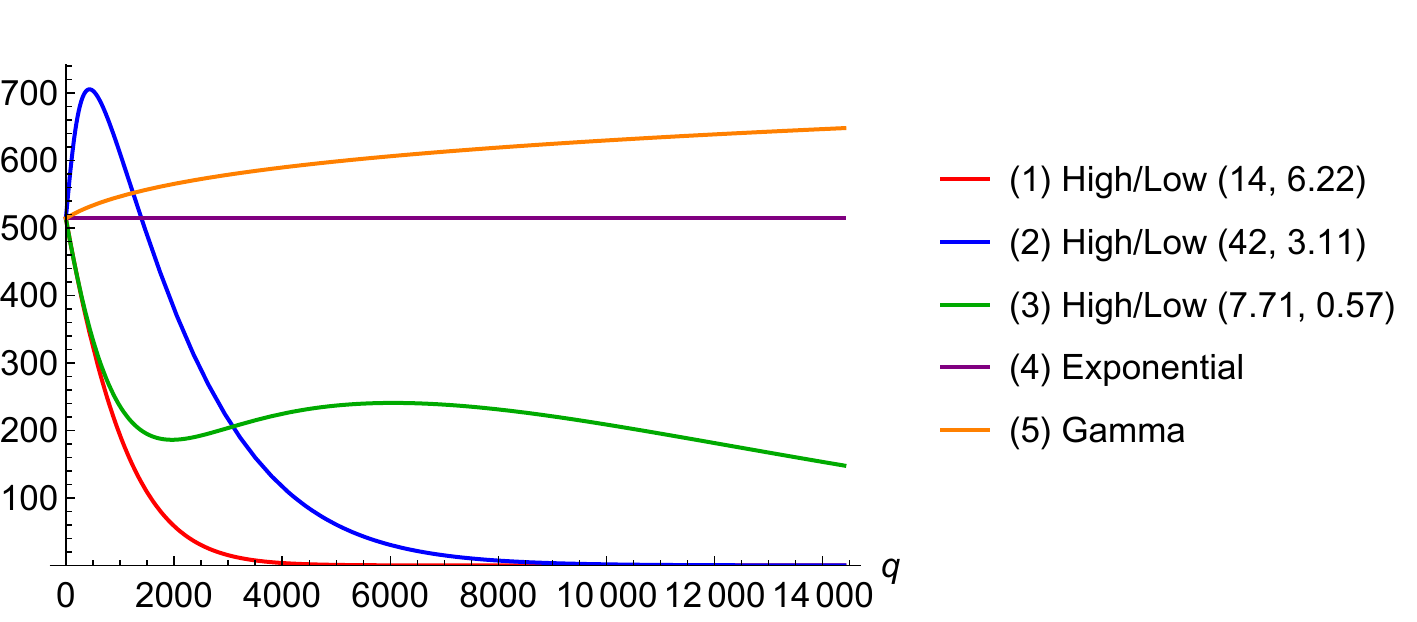}
\caption{Capacity of the minor road (veh/h) for several distributions of $T$.}
\label{fig:example4capacity}
\end{figure}

The previous example has illustrated that resampling, as described in B$_2$, has a positive impact on the capacity of the minor street. In this example we show that, under specific circumstances, this positive impact may be even bigger than expected. We show this by varying the probability distribution of the critical headway $T$, taking the following five distributions (all with $\mathbb{E}[T]=7$ seconds):
\begin{enumerate}
\item $T$ is equal to 14 seconds with probability $1/10$, and 6.22 seconds with probability $9/10$. This distribution, referred to as High/Low $(14, 6.22)$ in Figure\ \ref{fig:example4capacity}, is the same as in Example 1, but with more realistic values.
\item $T$ is equal to 42 with probability $1/10$, and 3.11 with probability $9/10$. This is a similar distribution as the previous, but with more extreme values.
\item $T$ is equal to 0.57 with probability $1/10$, and 7.71 with probability $9/10$. One out of ten samples is extremely small, instead of extremely large as in the previous distribution.
\item $T$ is exponentially distributed with parameter $1/7$.
\item $T$ has a gamma distribution with shape parameter $1/2$ and rate $1/14$.
\end{enumerate}

In Figure\ \ref{fig:example4capacity} we have plotted the capacity of the minor road as a function of $q$, the traffic intensity of the main road. Several conclusions can be drawn. First, we notice that in all High/Low distributions, the capacity drops to zero as $q$ increases towards infinity. However, when zooming in at $q$ close to zero, it turns out that models (1) and (3) have an immediate capacity drop, whereas model (2) actually \emph{increases} in capacity up to $q\approx 437$. A possible explanation for this striking phenomenon, that we also encountered in the previous numerical example, is that in the third model the critical headway is either extremely low, or extremely high. When a driver has a low critical capacity (which is quite likely to happen), he experiences no delay before leaving the minor street anyway. On the contrast, when a driver has a high critical capacity, he will have to resample, meaning that he has to wait for the next car to pass on the major street. If $q$ increases, the frequency of resampling increases, meaning that he has to wait less before getting a new chance to obtain a small critical headway.

Model (3) stands out because of its strange shape: first it drops to a local minimum around $q\approx1965$. After attaining this minimum, the capacity increases until $q\approx 6055$. Beyond this value, the capacity slowly drops to zero.

Another extreme case is the exponential distribution (4), which, due to its memoryless property, has a \emph{constant} capacity, not depending on the traffic flow on the major street.

Clearly, the most paradoxical case is (5), the Gamma distribution with a shape parameter less than one. If $T$ has this particular distribution, it can be shown that the capacity of the minor road keeps on \emph{increasing} as $q$ increases. Although this case, admittedly, may not be a realistic one, it again stresses that various sorts of counterintuitive phenomena may arise, and that one has to be very careful when applying heuristic reasoning.

\section{Impatience}\label{sect:impatience}

The goal of this section is to make the model more realistic by incorporating driver's impatience. As evidenced by Abou-Henaidy {\it et al.} \cite{abouhenaidy94}, drivers tend to grow more impatient as the number of rejected gaps increases. This impatience may result in an increased willingness to accept smaller gaps. 
To the best of our knowledge, \cite{abhishekcomsnets2016} was the first to present new results for gap acceptance models that include impatience \emph{and} randomness in the critical headways. We included a brief recap of some of the main results from \cite{abhishekcomsnets2016} and conclude this section by showing an example where impatience has a counterintuitive effect on the capacity of the minor road.

As discussed in Section \ref{sect:model}, we let the critical headway depend on the number of failed attempts. Denote by $T_j$ the critical headway for the $j$-th attempt to enter the main road ($j=1,2,\dots$). For models B$_1$ and B$_3$, we assume that $T_1 \geq T_2\geq \dots \geq T_{\text{min}}$ (more details below). Due to the resampling in model B$_2$, we cannot make this assumption for this model, but we can assume that $T_i \geq_{st} T_{i+1}$, where $\geq_{st}$ is used to denote that $T_i$ is stochastically greater than $T_{i+1}$. For each of the models B$_1$, B$_2$, and B$_3$, we compute the LST and the expectation of the ``service time'' $Y$, which leads to the capacity  of the minor road (Equation \eqref{fourc}). Note that the LST can also be used to find the transforms of the distributions of the queue length and the delay on the minor road. We refer the reader to \cite{abhishekcomsnets2016} for more details.

\vb

\noindent
\textbf{\boldmath B$_1$ {(constant gap):}} 
In this case, with fixed $T_j$, we use the memoryless property to show that the probability that the $j$-th attempt is successful is equal to ${\mathbb P}(\tau_q>T_j) = e^{-qT_j}$, where $\tau_q$ is an
exponentially distributed random variable with mean $1/q.$ It readily follows that the Laplace transform of the `service time', ${\mathbb E} [e^{-sY}]$, is
\begin{align}
\mathbb{E}[e^{-sY}]&=\sum_{k=0}^\infty \left(\prod_{j=1}^k {\mathbb E} [e^{-s\tau_q}1_{\{\tau_q<T_j\}}]\right)\,{\mathbb E} [e^{-sT_{k+1}}1_{\{\tau_q\geqslant T_{k+1}\}}]\label{lstB1}\\
&=\sum\limits_{k=0}^{\infty}\left(\frac{q}{s+q}\right)^ke^{-(s+q)T_{k+1}}\prod\limits_{j=1}^{k}(1-e^{-(s+q)T_j}).
\label{lstB1final}
\end{align}

Differentiation with respect to $s$ gives the expected value,
\begin{align*}
\mathbb{E}[Y]=
\sum\limits_{k=0}^{\infty}e^{-qT_{k+1}}\left[ \frac{k}{q}+T_{k+1}-\sum\limits_{i=1}^{k}\frac{T_ie^{-qT_i}}{1-e^{-qT_i}}\right]\prod\limits_{j=1}^{k}(1-e^{-qT_j}).
\end{align*}

Substitution of $T_1=T_2=\dots=:T$ in these expressions leads, after considerable simplification, to the expressions in the case without driver impatience from Section \ref{sect:analysis}.

\vspace{2mm}

\noindent
\textbf{\boldmath B$_2$ (sampling   per attempt):}
Despite the fact that $T_j$ is random now, we can still use the memoryless property of the gaps between successive cars on the major road in combination with the independence of the $T_j$, to argue that expression \eqref{lstB1} is also valid for this model. Note that ${\mathbb E} [e^{-s\tau_q}1_{\{\tau_q<T_j\}}]$ and ${\mathbb E} [e^{-sT_{k+1}}1_{\{\tau_q\geqslant T_{k+1}\}}]$ are different, though, yielding a slightly different result,
\[
\mathbb{E}[e^{-sY}]=\sum\limits_{k=0}^{\infty}\left(\frac{q}{s+q}\right)^k\mathbb{E}[e^{-(s+q)T_{k+1}}]\prod\limits_{j=1}^{k}(1-\mathbb{E}[e^{-(s+q)T_j}]).
\]
The mean service time follows from differentiating this expression:
\[\mathbb{E}[Y]=\sum\limits_{k=0}^{\infty}\left[ \mathbb{E}[T_{k+1}e^{-qT_{k+1}}]+\mathbb{E}[e^{-qT_{k+1}}]\left(\frac{k}{q}-\sum\limits_{i=1}^{k}\frac{\mathbb{E}[T_ie^{-qT_i}]}{1-\mathbb{E}[e^{-qT_i}]}\right)\right]\times \prod\limits_{j=1}^{k}(1-\mathbb{E}[e^{-qT_j}]).\]


\noindent
\textbf{\boldmath B$_3$ (sampling per driver):}   In this variant, every car driver samples a random $T_1$ at his first attempt. This (random) value determines the complete sequence of $T_2, T_3, \dots$. To model this kind of behavior, we introduce a sequence of functions $h_j(\cdots)$, for $j=1, 2, \dots$, such that the critical headway $T_j$ is defined as $T_j:=h_j(T_1)$, where $h_1$ is the identity function. 

Since only the first critical headway is random, we obtain the LST of $Y$ by conditioning on the value of $T_1$, and using \eqref{lstB1final}. We obtain
\begin{equation*}
\mathbb{E}[e^{-sY}]=\mathbb{E}\left[\sum\limits_{k=0}^{\infty}\left(\frac{q}{s+q}\right)^ke^{-(s+q)T_{k+1}}\prod\limits_{j=1}^{k}(1-e^{-(s+q)T_j}),\right]
\end{equation*}
and mean service time
\[
\mathbb{E}[Y]=\mathbb{E}\left[
\sum\limits_{k=0}^{\infty}e^{-qT_{k+1}}\left( \frac{k}{q}+T_{k+1}-\sum\limits_{i=1}^{k}\frac{T_ie^{-qT_i}}{1-e^{-qT_i}}\right)\times\prod\limits_{j=1}^{k}(1-e^{-qT_j})\right].
\]

A natural question is whether the ordering of the capacities observed in the case without impatience, $\bar\lambda_2\geqslant \bar\lambda_1\geqslant \bar\lambda_3$, also holds in the case with impatience. Since the expressions for $\mathbb{E}[Y]$ in the various models, and, as a consequence, the expressions for $\bar\lambda_1, \bar\lambda_2$, and $\bar\lambda_3$ are more complicated and less insightful than in the previous section, we will do a numerical experiment below to study the ordering of the capacities if drivers grow impatient.

\subsection{Example 3: impatience}

We now revisit model (1) from Example 2, but introducing driver impatience into the model, in the following specific form:
\begin{equation}
T_{k+1}=\alpha (T_k-\Delta)+\Delta, \qquad k = 1, 2, \dots; 0<\alpha<1,\label{impatienceexample}
\end{equation}
which means that the critical headway decreases in every next attempt, approaching the limiting value of $\Delta$. The parameter $\alpha$ determines the speed at which the patience decreases. In scenario 1 all $T_k$ are fixed, with $T_1=7$ seconds. In scenario 2, each of the $T_k$ is a random variable, with $T_1$ equal to $6.22$ or $14$ seconds, with probability $9/10$ and $1/10$ respectively. The distribution of $T_k$ for $k>1$ can be determined from  \eqref{impatienceexample}. Note that the impatience is a new random sample at each attempt, independent of the value of $T_{k-1}$. Scenario 3, as before, is similar to scenario 2, but each driver samples a random impatience $T_1$ exactly once. The value of $T_1$ (which is again either $6.22$ or $14$ seconds) determines the whole sequence of critical gap times at the subsequent attempts according to \eqref{impatienceexample}.

\begin{figure}[!ht]
\parbox{0.49\linewidth}{\centering
\includegraphics[width=\linewidth]{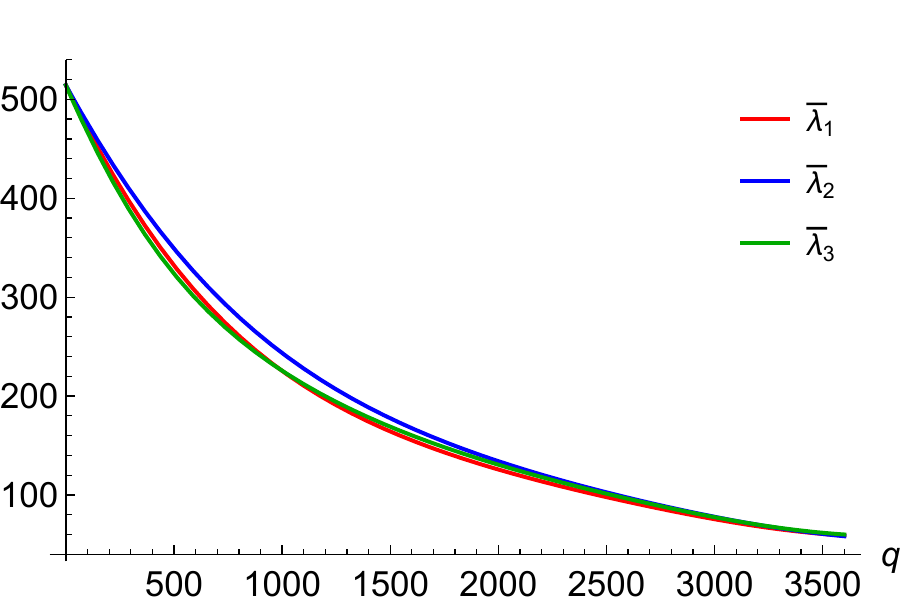}\\
\scriptsize (a) $\alpha=9/10, \Delta=4$ sec.}
\parbox{0.49\linewidth}{\centering
\includegraphics[width=\linewidth]{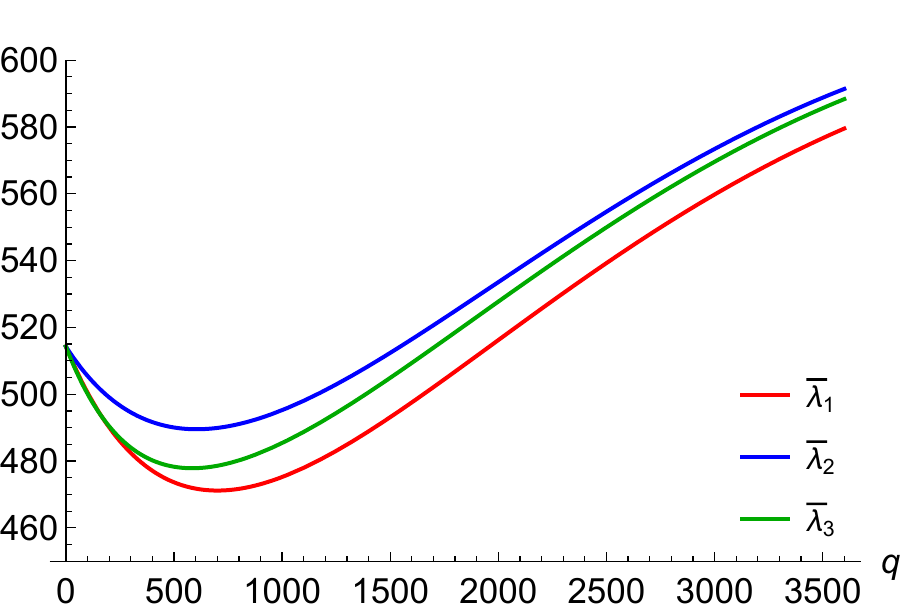}\\
\scriptsize (b) $\alpha=6/10, \Delta=1$ sec.}
\caption{Capacity of the minor street (veh/h) as a function of the flow rate on the main road (veh/h) in Example 3.}
\label{fig:example2capacity}
\end{figure}

Figure\ \ref{fig:example2capacity}(a) shows the capacity as a function of $q$, when $\alpha=9/10$ and $\Delta=4$ seconds. It is noteworthy that the strict ordering that was observed in the case without impatience, is no longer preserved, even though $\mathbb{E}[T_k]$ is the same in all scenarios, for $k$ fixed. Another interesting phenomenon, depicted in Figure\ \ref{fig:example2capacity}(b), occurs when we decrease the parameter values to $\alpha=6/10$ and $\Delta=1$ second. Now, the capacity actually \emph{increases} when $q$ exceeds a certain threshold. Due to the increase in $q$, gaps between cars on the major road will be smaller in general, but apparently the benefit of having a (much) lower critical headway at each attempt outweighs the disadvantage of having smaller gaps.

\section{Markov platooning on the major road}\label{sect:markovplatooning}

Assuming a Poisson arrival process on the major road is realistic in periods of free traffic flow, where it is assumed that any vehicle does not affect vehicles behind it. To make the model more realistic, Heidemann and Wegmann \cite{heidemann97}, relying on results by Tanner \cite{tanner62}, propose a general framework based on gap-block models.
In such models, vehicles form platoons which arrive according to a Poisson process.  The lengths of these platoons are i.i.d.\ random variables with a general distribution, which can be suitable chosen such that it matches real-life clustering behavior.
Wu \cite{wu2001} observed that, in practice, traffic flow in the major stream can have up to four different regimes: free space (no vehicles), free flow (single vehicles), bunched traffic (platoons of vehicles), and queueing. By conditioning on the current regime, he applies the  framework of \cite{heidemann97} to set up a heuristic argument that provides a more general capacity formula that is valid under all four regimes; we return to this approach below.

\begin{figure}[!ht]
\begin{center}
\includegraphics[width=\linewidth]{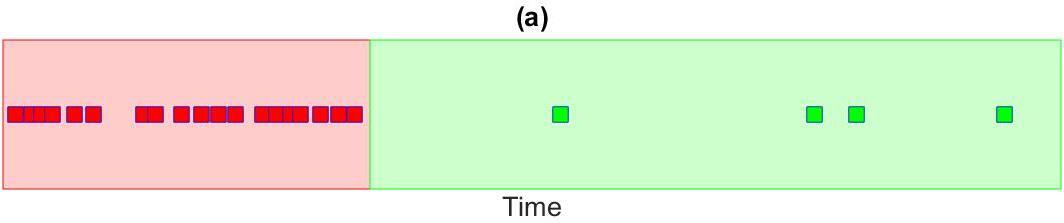}
\end{center}
\begin{center}
\includegraphics[width=\linewidth]{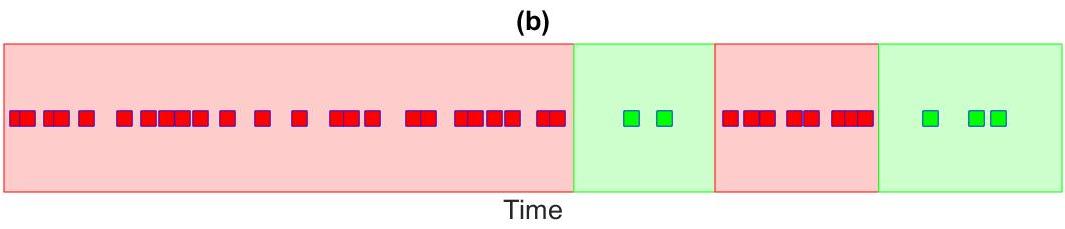}
\end{center}
\caption{Simulated examples of two MMPP's with two background states. On the horizontal axis we depict the time, while the squares (red or green) mark the arrivals.  In (a), we have chosen $\mu_1=1/20, \mu_2=1/40$ and arrival rates $q_1=1, q_2=1/15$ vehicles per time unit. In (b) we use $\mu_1=1/20, \mu_2=1/20$ and arrival rates $q_1=1, q_2=1/5$ vehicles per time unit. The red areas indicate that the background process is in state 1 (more platooning) and the green areas correspond to state 2 (less platooning).}
\label{fig:intersectionarrivals}
\end{figure}

In this section we introduce a novel way to model different traffic-flow regimes on the major road, using a well-established method to model dependence between successive inter-arrival times. We assume that the arrival process on the major road is modeled by a {\it Markov modulated Poisson process} (MMPP).
In an MMPP arrivals are generated at a Poisson rate $q_i$ when an exogeneous, autonomously evolving continuous-time Markov process (commonly referred to as the {\it background process}) is in state $i$. We denote by  $d \in\{1,2,\dots\}$ the number of states of the background process (where $d=1$ corresponds to a non-modulated, ordinary Poisson process). We assume the background process to be irreducible; the corresponding stationary distribution is given by the vector $\pi$.
In the sequel we denote by $M=(\mu_{ij})_{i,j=1}^d$ the transition rate matrix of background process, and define $\mu_i:=-\mu_{ii}$. Therefore, an MMPP allows different traffic-flow regimes on the major road.
For example, in Figure \ref{fig:intersectionarrivals}, we show the arrival patterns of two MMPP's, each with two background states. The red squares mark arrivals during the high traffic intensity ($q_1$), while the green squares mark arrivals during the low intensity ($q_2$). It can be seen that platoons are generally longer when the background process is in state $1$, corresponding to a high arrival rate. Additionally, we observe in Figure \ref{fig:intersectionarrivals}(a) that the background process stays longer in state $2$ ($\mu_2=1/40$) than in state $1$ ($\mu_1=1/20$).
Another difference between the two sub-figures is that we choose $q_2=1/15$ in Figure \ref{fig:intersectionarrivals}(a) and $q_2=1/5$ in Figure \ref{fig:intersectionarrivals}(b). This explains why, in state 2, we see no platooning at all in Figure  \ref{fig:intersectionarrivals}(a), but Figure \ref{fig:intersectionarrivals}(b) still shows some mild platoon forming.

The main objective of this section is to develop methods that determine the capacity of the minor road under MMPP arrivals on the major road, for the models B$_1$ up to B$_3$. Because of this focus on the capacity, we can simplify the model by taking away the queueing aspect on the minor road, assuming that this road is {\it saturated}: there are \emph{always} low-priority cars waiting for gaps. The reason underlying this reduction is that capacity is a quantity that corresponds to stability of the associated queue, and stability essentially amounts to the queue being able to process all input in the long run.

The capacity, to be denoted by $\bar\lambda$, is the ratio of the mean number of arrived cars in a cycle (which we define below) to the mean duration of a cycle, which equals (due to renewal theory) the number of cars that can be served per unit time. The system is stable when $\lambda$, the arrival intensity on the minor road, is less than $\bar\lambda$. Again, we distinguish between the three behavior types B$_1$, B$_2$, and B$_3$ introduced in Section \ref{sect:model}, each with its own capacity $\bar\lambda_i$, for $i=1,2,3$.

\vb

Our objective is to assess the impact of the three types of the driver's behavior on stability. As we have seen before, the capacity can be interpreted as the reciprocal of the time it takes for an arbitrary car to cross the major road. At first sight, the following procedure seems to provide us with $\bar\lambda$. Define $S_i$ as the time it takes for an arbitrary car to cross the major road, given the background process is in state $i$ when the car (which has reached the head of the queue) starts his attempt. Recalling that $\pi_i$ represents the long-run fraction of time that the background process resides in state $i$, it is tempting to conclude that the capacity would equal
\begin{equation}\label{eqn:wrongcapacityformula_1}
\sum_{i=1}^d \frac{\pi_i}{ {\mathbb E}[S_i]}.
\end{equation}
Alternatively, one might try to first take a weighted average of the mean service times, and then take the reciprocal to find the capacity,
\begin{equation}\label{eqn:wrongcapacityformula_2}
 \frac{1}{\sum_{i=1}^d \pi_i{\mathbb E}[S_i]}.
\end{equation}
There is, however, a conceptual mistake in these (na\"{\i}ve) approaches. It is true that $\pi$ is indeed the distribution of the background process that is seen by cars that arrive at the queue, due to the well-known PASTA property. The distribution seen by the car that has reached the head of the queue, however,  differs from $\pi$. To see this, think of the extreme case in which $q_1=q>0$ has some moderate value, and $q_2=M$ is large. For large values of $M$, only cars who find the queue empty, may start their attempt while the background process is in state $2$.

This reasoning illustrates how careful one should be when weighing capacities that belong to different regimes by the fractions of time in which those regimes apply. A very similar decomposition approach was  followed by Wu \cite{wu2001}; he distinguishes four different regimes, as described above, each with an own capacity, and those are combined into a single capacity. The formulas obtained by Wu \cite{wu2001} likely provide a reasonable indication of the capacity across a wide range of parameters, but there are also many cases in which the approach fails to do so. Later on, we provide an example which illustrates what errors may result when following
the na\"{\i}ve approaches.

\vb

\noindent
\textbf{\boldmath B$_1$ {(constant gap):}}
In this model, every driver on the minor road needs the same constant critical headway $T$ to enter the major road.  In our analysis we use the renewal reward theorem, which entails that the capacity can be written as the mean number of cars arriving in a regenerative  cycle  divided by the mean duration of that cycle.
For our purposes, an appropriate definition of a cycle is: the time elapsed between two consecutive epochs such that (i) the background process is in a reference state (say state $1$), and (ii) a service is completed (i.e., a low-priority car is served).

To make our model Markovian,  we approximate this deterministic $T$ by an Erlang random variable with $k$ phases of average length $T/k$. It is well known that a deterministic $T$ can be approximated by the sum of $k$ independent exponential random variables, each with parameter $\kappa:=k/T$, with $k$ large; to see this, observe that this Erlang random variable has mean $T$ (as desired), and variance
${k}/{\kappa^2} = {T^2}/{k},$
which goes to $0$ as $k$ grows large. In the sequel we write $\varrho_i:=\mu_i+q_i+\kappa$. Define $h_{ij}$ as the mean number of cars that is served till the cycle ends, given that the current state of the background process is $i\in\{1,\ldots,d\}$ and the car in  service has finished $j\in\{0,\ldots,k-1\}$ phases of the Erlang distribution. To find the  mean number of arrived cars in a cycle, we need to find $h_{10}$. This can be done as follows.

Relying on `standard Markovian reasoning', by conditioning on the first jump,
\[h_{1j}  =\sum_{\ell\not= 1}\frac{\mu_{1\ell}}{\varrho_1} h_{\ell j}+
\frac{q_1}{\varrho_1} h_{10} +
\frac{\kappa}{\varrho_1} \left(h_{1,j+1} 1_{\{j<k-1\}} + 1_{\{j=k-1\}} \right)
.\]
In addition, for $i\not=1$,
\[h_{ij}  =\sum_{\ell\not= i}\frac{\mu_{i\ell}}{\varrho_i} h_{\ell j}+
\frac{q_i}{\varrho_i} h_{i0} +
\frac{\kappa}{\varrho_i} \left(h_{i,j+1} 1_{\{j<k-1\}} + (1+h_{i0})1_{\{j=k-1\}} \right)
.\]
This can be written as a linear system of $dk$ equations with $dk$ unknowns of the form $A\vec{h} = \vec{b}$, where entries of the matrix $A=[a_{mn}]$, $\vec{h}$ and $\vec{b}=[b_{m}]$ are given as follows, with $i=\lceil m/k\rceil$,
\[ a_{mn} = \left\{
\begin{array}{l l}
-{\displaystyle \frac{\kappa}{{\varrho_{i}}}}, & \quad \text{if $n=m+1$ and $m\neq k,2k,...,dk$;}\\
&\vspace{-4mm}\\
-{\displaystyle \frac{q_{i}}{{\varrho_{i}}}}, & \quad \text{if $n=(i-1)k+1$ and $m\neq 1,k+1,2k,2k+1,3k,...,(d-1)k+1,dk$;}\\
&\vspace{-4mm}\\
-{\displaystyle \frac{\kappa + q_{i}}{{\varrho_{i}}}}, & \quad \text{if $n=(i-1)k+1$ and $m=2k,...,dk$;}\\
&\vspace{-4mm}\\
-{\displaystyle \frac{\mu_{i,\ell +1}}{{\varrho_{i}}}}, & \quad \text{if $n=(\ell -i+1)k+m$ and $\ell \in \{0,1,...,d-1\}\setminus\{i-1\}$;}\\
&\vspace{-4mm}\\
1-{\displaystyle \frac{q_{i}}{{\varrho_{i}}}}, & \quad \text{if $n=m$ and $m=1,k+1,...,(d-1)k+1$;}\\
1, & \quad \text{if $n=m$ and $m\neq 1,k+1,...,(d-1)k+1$;}\\
0, & \quad \text{else,}
\end{array} \right.\]
$\vec{h}=[h_{10},h_{11},...,h_{1 ,k-1},h_{20},h_{21},...,h_{2, k-1},...,h_{d0},h_{d1},...,h_{d ,k-1}]^{\rm T}$, and
\[ b_{m} = \left\{
\begin{array}{l l}
{\displaystyle \frac{\kappa}{{\varrho_{\ell }}}}, & \quad \text{if $m=\ell k, \ell =1,2,...,d$ }\\
0, & \quad \text{else.}
\end{array} \right.\]

It is noted that $|a_{mm}|=\sum_{n\neq m}|a_{mn}|$ for all $m\neq k$ and for $m=k$, $|a_{kk}|>\sum_{n\neq k}|a_{kn}|$. Therefore, the matrix $A$ is  weak diagonally dominant with one row is strictly dominant and $A$ is also irreducible matrix, and hence $A$ is invertible \cite{horn}. Therefore, the solution of the system of equations $A\vec{h} = \vec{b}$ is $\vec{h} = A^{-1}\vec{b}$. We thus find the desired quantity  $h_{10}$.

\vb

To determine the capacity we need, in addition to the  mean number of arrived cars in a cycle, also the mean duration of a cycle. To this end
we define $\tau_{ij}$ as the mean time till the end of the current cycle, given that the current state of the background process is $i\in\{1,\ldots,d\}$ and the car in service has finished $j\in\{0,\ldots,k-1\}$ phases of the Erlang distribution. The objective is now to find the mean duration of a cycle, which is given by $\tau_{10}$.

Similarly to the procedure we set up above,
\[\tau_{1j}  =\frac{1}{\varrho_1}+
\sum_{\ell\not= 1}\frac{\mu_{1\ell}}{\varrho_1} \tau_{\ell j}+
\frac{q_1}{\varrho_1} \tau_{10} +
\frac{\kappa}{\varrho_1} \tau_{1,j+1} 1_{\{j<k-1\}}
.\]
In addition, for $i\not=1$,
\[\tau_{ij}  =\frac{1}{\varrho_i} +
\sum_{\ell\not= i}\frac{\mu_{i\ell}}{\varrho_i} \tau_{\ell j}+
\frac{q_i}{\varrho_i} \tau_{i0} +
\frac{\kappa}{\varrho_i} \left(\tau_{i,j+1} 1_{\{j<k-1\}} + \tau_{i0}1_{\{j=k-1\}} \right)
.\]
Also this system can be written as $dk$ linear equations with $dk$ unknowns. More precisely, with $\vec{\tau}=[\tau_{10},\tau_{11},...,\tau_{1, k-1},\tau_{20},\tau_{21},...,\tau_{2, k-1},...,\tau_{d0},\tau_{d1},...,\tau_{d ,k-1}]^{\rm T}$, we have
$A\vec{\tau} = \vec{c}$ with $A$ as defined before and  \[c_{m} =
\frac{1}{\varrho_{\ell }}\:\:\:\mbox{for $(\ell -1)k+1 \leq m\leq \ell k, \:\:\mbox{and}\:\: \ell =1,2,...,d.$}\]
We already proved that $A$ is invertible, and therefore the unique solution of the system of equations $A\vec{\tau} = \vec{c}$ is $\vec{\tau} =A^{-1} \vec{c}$. We thus find $\tau_{10}$.

The capacity of this system can now be evaluated as $\bar\lambda_1:=h_{10}/\tau_{10}$, meaning that the stability condition of the low-priority queue is $\lambda< \bar\lambda_1$. In the numerical procedure, the value of  $k$ should be chosen large, to ensure that the Erlang distribution is sufficiently `close-to-deterministic'.

\vb

\noindent
\textbf{\boldmath B$_2$ (sampling   per attempt):}
Appointed out before, in this behavior type every driver samples a `fresh' random $T$ for every attempt to enter  the major road. Let us assume that the gap size $T$ equals some deterministic $T_n$ with probability $p_n$; below we present the computational procedure for $n\in\{1,2\}$, but it can be extended in an evident manner to the situation in which $T$ can attain more than $2$ possible values.  Analogously to what we did in the procedure to evaluate the capacity for B$_1$, we approximate $T_n$ by an Erlang random variable with $k_n$ phases; each of the phases is exponentially distributed with parameter $\kappa_n=k_n/T_n$.

We  write $\varrho^{(n)}_i:=\mu_i+q_i+\kappa_n.$ Let $h_{ij}^{(n)}$ be the mean number of cars that is served till the cycle  ends, given that the current state of the background process is $i\in\{1,\ldots,d\}$, the car in service has gap size $T_n$ and the car in service has finished $j\in\{0,\ldots,k_n-1\}$ phases. We wish to find $h_{10}$ where
\begin{equation}\label{hi0}
h_{i0}=p_1h_{i0}^{(1)}+p_2h_{i0}^{(2)}$ for $i=1,2.\end{equation}
Then
\[h^{(n)}_{1j}  =\sum_{\ell\not= 1}\frac{\mu_{1\ell}}{\varrho^{(n)}_1} h^{(n)}_{\ell j}+
\frac{q_1}{\varrho^{(n)}_1} h_{10} +
\frac{\kappa_{n}}{\varrho^{(n)}_1} \left(h^{(n)}_{1,j+1} 1_{\{j<k_n-1\}} + 1_{\{j=k_n-1\}} \right)
.\]
Observe how the resampling is incorporated in this system: when an attempt has failed a `fresh' new gap size is sampled, explaining the $h_{10}$ (rather than $h^{(n)}_{10}$) in the right hand side.
In addition, for $i\not=1$,
\[h^{(n)}_{ij}  =\sum_{\ell\not= i}\frac{\mu_{i\ell}}{\varrho^{(n)}_i} h^{(n)}_{\ell j}+
\frac{q_i}{\varrho^{(n)}_i} h_{i0}  +
\frac{\kappa_{n}}{\varrho^{(n)}_i} \left(h^{(n)}_{i,j+1} 1_{\{j<k_n-1\}} + (1+h_{i0})1_{\{j=k_n-1\}} \right)
.\]
The $h_{i0}$ in the right hand side of the previous display corresponds with the event that an attempt has succeeded, after which   a  new gap size is sampled.

The above equations can be written as a linear system of the type  $A\vec{h} = \vec{b}$ for a matrix $A$ and vector $\vec{b}$  (which evidently differ from the matrix $A$ and vector $\vec{b}$ that were used in the model
B$_1$) consisting of $d(k_1+k_2)$ equations with $d(k_1+k_2)$ unknowns.
With the same argument as we have used for  B$_1$, it follows that the coefficient matrix $A$ is invertible. Using (\ref{hi0}), this facilitates the computation of $\vec{h}$ and in particular the desired quantity $h_{10}$ (from  $h_{10}=p_1 {h}^{(1)}_{10}+p_2 {h}^{(2)}_{10}$).

\vb

We then define $\tau^{(n)}_{ij}$ as the mean time till the current cycle ends, given that the current state of the background process is $i\in\{1,\ldots,d\}$, the car in the service has gap size $T_n$ and the car in service has finished $j\in\{0,\ldots,k_n-1\}$ phases. The objective is to set up a numerical procedure to evaluate $\tau_{10}$ where $\tau_{i0}=p_1\tau_{i0}^{(1)}+p_2\tau_{i0}^{(2)}$ for $i=1,2$. Using the same argumentation as above,
\[\tau^{(n)}_{1j}  =\frac{1}{\varrho^{(n)}_1}+\sum_{\ell\not= 1}\frac{\mu_{1\ell}}{\varrho^{(n)}_1} \tau^{(n)}_{\ell j}+
\frac{q_1}{\varrho^{(n)}_1} \tau_{10}  +
\frac{\kappa_{n}}{\varrho^{(n)}_1} \tau^{(n)}_{1,j+1} 1_{\{j<k_n-1\}}.\]
In addition, for $i\not=1$,
\[\tau^{(n)}_{ij}  =\frac{1}{\varrho^{(n)}_i}+\sum_{\ell\not= i}\frac{\mu_{i\ell}}{\varrho^{(n)}_i} \tau^{(n)}_{\ell j}+
\frac{q_i}{\varrho^{(n)}_i} \tau_{i0} +
\frac{\kappa_{n}}{\varrho^{(n)}_i} \left(\tau^{(n)}_{i,j+1} 1_{\{j<k_n-1\}} + \tau_{i0}1_{\{j=k_n-1\}} \right)
,\]
with
$
\tau_{i0}=p_1\tau_{i0}^{(1)}+p_2\tau_{i0}^{(2)}$ for $i=1,2.$
Again, this system can be written as a linear system of $d(k_1+k_2)$ equations with $d(k_1+k_2)$ unknowns, say $A\vec{\tau} = \vec{c}$, with $A$ as above (and hence invertible). Therefore, the solution of the system of equations $A\vec{\tau} = \vec{c}$ is $\vec{\tau} =A^{-1} \vec{c}$, and we can compute  $\tau_{10}=p_1{\tau}^{(1)}_{10}+p_2{\tau}^{(2)}_{10}$.
The capacity of the low-priority queue under B$_2$ is therefore  $\bar\lambda_2={h_{10}}/{\tau_{10}}$.

\vb

\vspace{2mm}

\noindent
\textbf{\boldmath B$_3$ (sampling   per driver):} We finally consider the model with consistent behavior, i.e., each driver sticks to the gap size he or she initially sampled.  The procedure is similar to the ones we developed for B$_1$ and B$_2$, and therefore we restrict ourselves to the main steps.

Define, as before,
$h_{i0}=p_1h_{i0}^{(1)}+p_2h_{i0}^{(2)}$ for $i=1,2.$
The mean number of cars served during the cycle follows from
\[h^{(n)}_{1j}  =\sum_{\ell\not= 1}\frac{\mu_{1\ell}}{\varrho^{(n)}_1} h^{(n)}_{\ell j}+
\frac{q_1}{\varrho^{(n)}_1} h^{(n)}_{10} +
\frac{\kappa_{n}}{\varrho^{(n)}_1} \left(h^{(n)}_{1,j+1} 1_{\{j<k_n-1\}} + 1_{\{j=k_n-1\}} \right)
;\]
it is instructive to compare this equation with the corresponding one for B$_2$: when the attempt has failed the gap size is {\it not} resampled.
Also, for $i\not=1$, along the same lines,
\[h^{(n)}_{ij}  =\sum_{\ell\not= i}\frac{\mu_{i\ell}}{\varrho^{(n)}_i} h^{(n)}_{\ell j}+
\frac{q_i}{\varrho^{(n)}_i} h^{(n)}_{i0} +
\frac{\kappa_{n}}{\varrho^{(n)}_i} \left(h^{(n)}_{i,j+1} 1_{\{j<k_n-1\}} + (1+h_{i0})1_{\{j=k_n-1\}} \right)
;\]
resampling is only done when an attempt has been successfully completed.

Similarly, the system of equations for the mean cycle length is
\[\tau^{(n)}_{1j}  =\frac{1}{\varrho^{(n)}_1}+\sum_{\ell\not= 1}\frac{\mu_{1\ell}}{\varrho^{(n)}_1} \tau^{(n)}_{\ell j}+
\frac{q_1}{\varrho^{(n)}_1} \tau^{(n)}_{10} +
\frac{\kappa_{n}}{\varrho^{(n)}_1} \tau^{(n)}_{1,j+1} 1_{\{j<k_n-1\}},
.\]
and, for $i\not=1$,
\[\tau^{(n)}_{ij}  =\frac{1}{\varrho^{(n)}_i}+\sum_{\ell\not= i}\frac{\mu_{i\ell}}{\varrho^{(n)}_i} \tau^{(n)}_{\ell j}+
\frac{q_i}{\varrho^{(n)}_i} \tau^{(n)}_{i0} +
\frac{\kappa_{n}}{\varrho^{(n)}_i} \left(\tau^{(n)}_{i,j+1} 1_{\{j<k_n-1\}} + \tau_{i0}1_{\{j=k_n-1\}} \right)
,\]
with
$
\tau_{i0}=p_1\tau_{i0}^{(1)}+p_2\tau_{i0}^{(2)}$ for $i=1,2.$
The linear system can solved as before, yielding $h_{10}$ and $\tau_{10}$.
Therefore, the capacity of the system can be evaluated as $\bar\lambda_3={h_{10}}/{\tau_{10}}$.

\vb

\subsection{Example 4: the impact of Markov platooning}

The purpose of this collection of numerical examples is to exhibit specific, interesting features of gap acceptance models that relate to the impact of Markov platooning. In the literature it has already been observed that platoon forming on the major road may have a positive impact on the capacity of the minor road. For the first example, which is similar to Example 1 but now with Markov platooning, we compare the capacity of the minor road for the three behavior types B$_1$, B$_2$, and B$_3$. For the last two behavior types, we assume that a driver requires either a short gap of $T_1=3$ seconds, or an extremely long gap of $T_2=60$ seconds. Obviously these values are not chosen with the intention to mimic realistic behavior, but to point out extreme situations that might occur. For behavior type B$_1$, we take $T=p_1 T_1+p_2 T_2$ seconds long, where $p_2:=1-p_1$.

For these settings, we compare the model with and without Markov platooning. With platooning, we take $\mu_1=1/60$ and $\mu_2=1/240$, resulting in exponential periods of, on average, one minute  where the arrival rate on the major road is $q_1$, followed by exponential periods of, on average, four minutes, with arrival rate $q_2$. We assume a fixed ratio of $q_1$ and $q_2$, namely $q_1=3q_2$. The long-term average arrival rate equals \[\bar{q}:=\frac{q_1/\mu_1+q_2/\mu_2}{1/\mu_1+1/\mu_2}=
\frac{q_1\mu_2+q_2\mu_1}{\mu_1+\mu_2}.\] We compare the capacities with those obtained from the model without platooning, where we assume Poisson arrivals with rate $\bar{q}$.

\begin{figure}[!ht]
\parbox{0.49\linewidth}{\centering
\includegraphics[width=\linewidth]{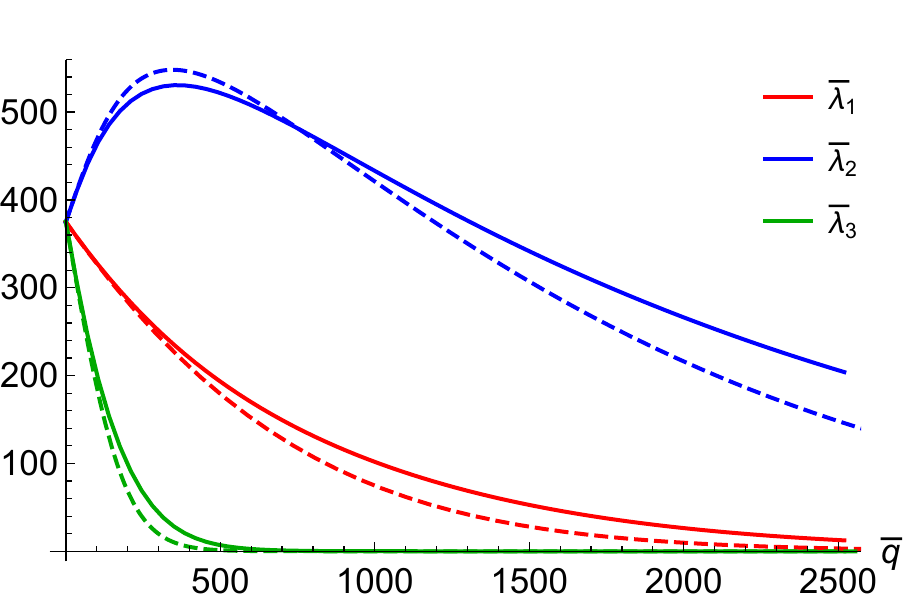}\\
\scriptsize (a) $p_1=0.9$}
\parbox{0.49\linewidth}{\centering
\includegraphics[width=\linewidth]{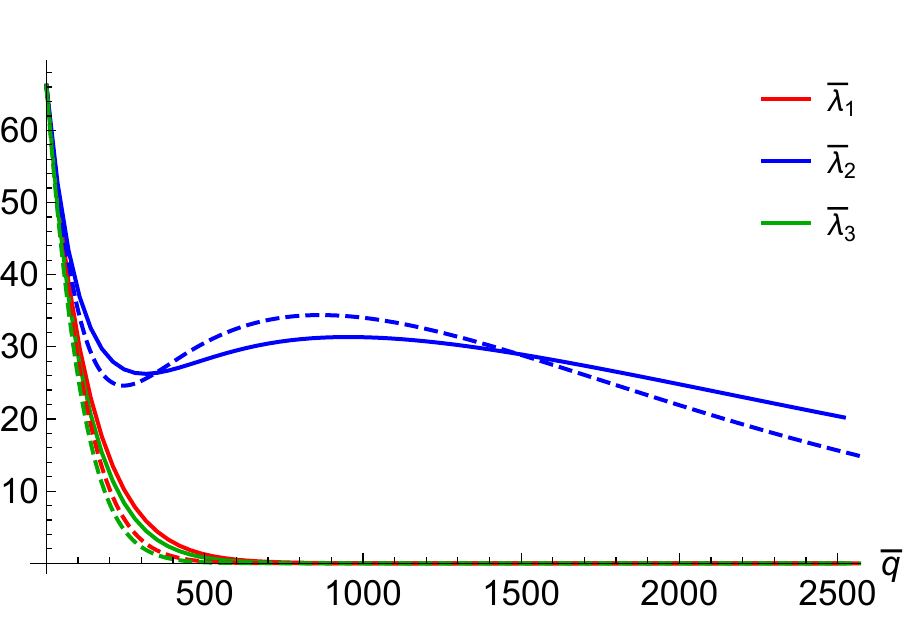}\\
\scriptsize (b) $p_1=0.1$.}
\caption{Capacity of the minor street (veh/h) as a function of the average flow rate on the main road (veh/h) in Example 4. The solid lines correspond to the model with Markov platooning; the dashed lines correspond to the model without platooning.}
\label{fig:ex5-1}
\end{figure}

Figure \ref{fig:ex5-1} depicts the capacity (veh/h) of the minor street as a function of $\bar{q}$, the average flow rate on the main road (veh/h), for $p_1=0.9$ and  $p_1=0.1$ respectively. As in the non-modulated case, we observe  the relation $\bar\lambda_2\geqslant\bar\lambda_1\geqslant\bar\lambda_3$. Due to the lack of explicit expressions for $\bar{\lambda}_1, \bar{\lambda}_2$, and $\bar{\lambda}_3$, we cannot prove the strict ordering now. We did, however, observe it in all numerical examples that we conducted, and conjecture the ordering to hold true in general.

Based on the results of this example (and many other examples that are not discussed in the present paper) we are inclined to believe that platooning has a positive effect on the capacity of the minor road, but \emph{only} for models B$_1$ and B$_3$. In a model with inconsistent behavior, it really depends on the model parameters whether platooning increases or decreases the capacity. This is nicely illustrated in Figure \ref{fig:ex5-1}(a) and even better in Figure~\ref{fig:ex5-1}(b).

\subsection{Example 5: platoon lengths}

In this example we fix the overall arrival rate on the major road, but we vary the platoon sizes. In more detail, we assume that $q_1=600$ veh/h and $q_2=2400$ veh/h. This means that phase 1 can be considered as a situation of moderate traffic (every 6 seconds a car passes), whereas phase 2 can be considered as one big platoon (on average every 1.5 seconds a car passes). The overall arrival rate $\bar{q}$ is fixed at 900 vehicles per hour, which implies that $\mu_1/\mu_2=1/5$. By varying the mean platoon length $1/\mu_2$ (in seconds) between 0 and 10, we will get better insight in the relation between platoon lengths and the capacity. Wegmann \cite[Section 5]{wegmann1991}  conducted a very similar experiment, varying the mean number of vehicles per bunch. He observed that the capacity increases with increasing variance of gaps.

\begin{figure}[!ht]
\parbox{0.49\linewidth}{\centering
\includegraphics[width=\linewidth]{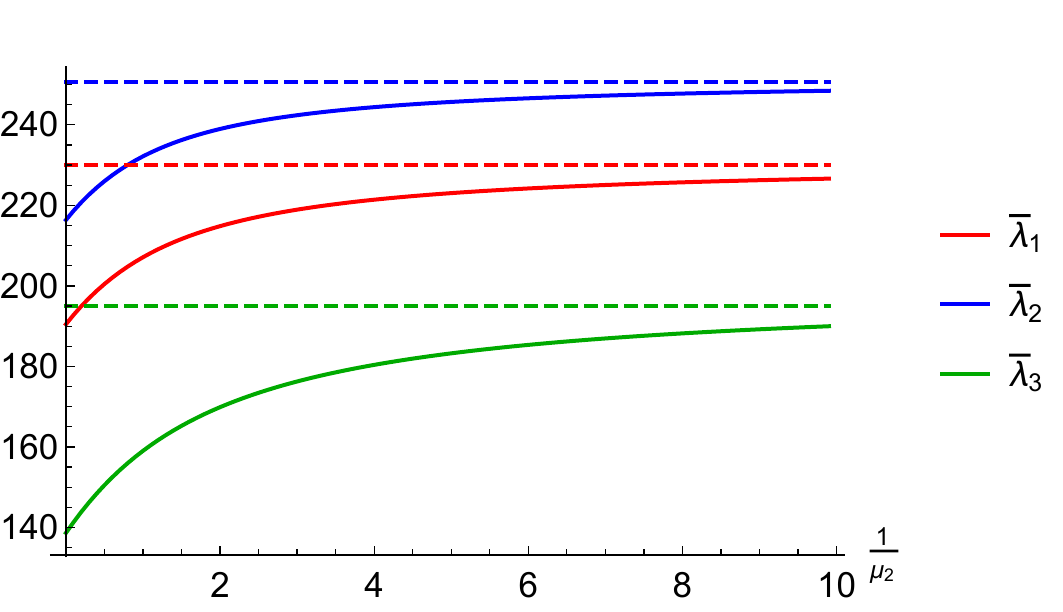}\\
\scriptsize (a) $T_1=6.22, T_2=14, p_1=0.9$}
\parbox{0.49\linewidth}{\centering
\includegraphics[width=\linewidth]{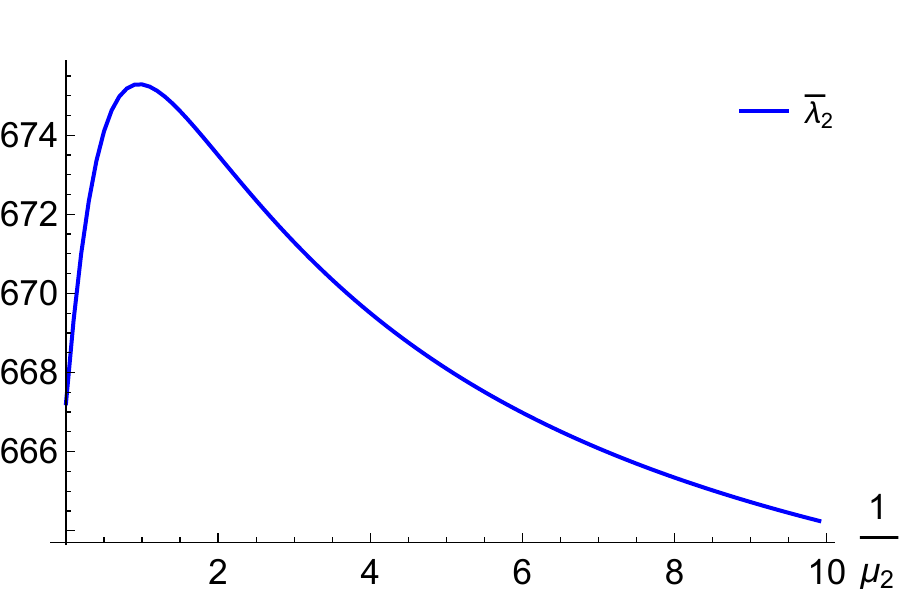}\\
\scriptsize (b) $T_1=3, T_2=60, p_1=0.9$}
\caption{Capacity of the minor street (veh/h) as a function of the mean platoon length (sec) in Example 5. The dashed lines in (a) indicate the limiting capacities for $\mu_2\downarrow 0$ while keeping the ratio $\mu_1/\mu_2$ fixed.}
\label{fig:ex6-1}
\end{figure}

We consider two different distributions for the critical headways. First, we consider the situation with $T_1=6.22$, $T_2=14$, and $p_1=0.9$, which can be considered as a quite realistic situation that we have used before. In Figure \ref{fig:ex6-1}(a) we show the results for behavior types B$_1$, B$_2$, and B$_3$. The relation between the capacity and the mean platoon length is in line with \cite[Figure 3]{wegmann1991}. Our numerical experiments confirm that this is indeed typical behavior for B$_1$, B$_2$, and B$_3$. Nevertheless, we want to show that it is possible to create a situation where model B$_2$ exhibits completely different behavior. When changing the distribution of the critical headway such that $T_1=3$ and $T_2=60$, we no longer see a monotonous relation between the capacity and the mean platoon length; see Figure \ref{fig:ex6-1}(b). Considering the fact that this inconsistent behavior type in combination with the extreme values for $T_1$ and $T_2$ might not be all too realistic, we do not find it likely that this type of behavior occurs in practical situations, but  the model shows that it is not entirely impossible. For completeness, we want to mention that under extreme circumstances such as mean platoon lengths of 1000 seconds, the capacity with consistent behavior B$_3$ will also exhibit a drop, but not as drastically as in Figure \ref{fig:ex6-1}(b).

The final conclusion that can be drawn from this example, is that one should be cautious when developing capacity estimates based on Equations \eqref{eqn:wrongcapacityformula_1} and \eqref{eqn:wrongcapacityformula_2}. This type of reasoning may create a substantial bias, due to the fact that the vehicle at the head of the queue typically does not see the background process in equilibrium. It is noted that such argumentation  underlies the capacity formulae in e.g.\ \cite{wu2001}, where
the capacity is calculated  by conditioning on the state of the background process, i.e., the state of the traffic on the major road (free space, free flow, bunching, or queueing). This example, and also Wegmann's example, clearly show  that there is a clear dependency between the mean platoon size and the capacity. The parameters in these examples are carefully chosen, such that the steady-state distribution of the background process (the vector $\pi$) remains unchanged. In our case, the major road is in state ``free flow'' for a fraction $\pi_1=5/6$ of the time, and in state ``bunched'' for a fraction $\pi_2=1/6$ of the time. If one would use the na\"{\i}ve approach and determine ${\mathbb E}[S_1]$ and ${\mathbb E}[S_2]$ by considering two separate models with regular Poisson arrivals, with intensities respectively $q_1$ and $q_2$, and use Equation \eqref{eqn:wrongcapacityformula_1}, the capacities for models B$_1$, B$_2$, and B$_3$, respectively, would be
\[
\bar{\lambda}_1=229.91, \bar{\lambda}_2 = 250.65, \bar{\lambda}_3=194.89,
\]
independent of $\mu_1$ and $\mu_2$. From Figure \ref{fig:ex6-1}(a) and Figure 3 in \cite{wegmann1991}, it is clearly visible that these values (indicated by the dashed lines in Figure \ref{fig:ex6-1}(a)) may differ substantially from the actual capacities. In fact, the capacities calculated from \eqref{eqn:wrongcapacityformula_1} can be interpreted as the limiting capacities from our MMPP model when $\mu_2\downarrow 0$ while keeping the ratio $\mu_1/\mu_2$ fixed. When using \eqref{eqn:wrongcapacityformula_2} to compute the capacities, one would obtain
\[
\bar{\lambda}_1=96.28, \bar{\lambda}_2 = 130.74, \bar{\lambda}_3=11.63,
\]
leading to even greater errors.

\section{Concluding remarks}\label{sect:conclusions}

Our main target in this work has been to investigate the impact of randomness in the critical headway on the capacity for traffic flows of low priority at a road intersection. For that, we have analyzed three versions of a queueing-theoretic model, each with its own dynamics. Special attention was paid to drivers' impatience under congested circumstances: the value of the critical headway decreases with subsequent attempts to cross the main road. In addition, we have provided a framework that allows the systematic evaluation of the effect of platooning on the primary road.

In our first model (B$_1$) we have assumed that the sequence of critical headways is a deterministically decreasing sequence, and that all cars use same sequence. In the second model (B$_2$), we let each car sample new values for the critical headway, according to a {\em stochastically decreasing} sequence. In the third model (B$_3$) we sample the first value for the critical headway for each car, but then use a deterministic decreasing sequence throughout the attempts of a particular car.

Unlike studies that appeared before, we focus on assessing the impact of the drivers' behaviors on the capacity for the low priority flow. Our main observation is that randomness has a strong impact on the capacity. More specifically, the capacity region depends on the entire distribution(s) of the critical headway durations, and not only of the mean value(s). We also observe that resampling of the critical headway values has a benign impact on the capacity of the minor road. To make our model more realistic than existing one, we have included two additional features into the model: impatience (the longer the driver has to wait, the lower the critical headway) and platooning (modelling the fluctuations in the traffic density at the primary road). Our results show that various counterintuitive and sometimes paradoxical phenomena may occur, and that heuristically developed guidelines should be handled with care.

We are currently investigating several extensions to our model. The analysis of our models carries over to variations of it in which, for example, cars may not need the entire duration of their critical headway to cross the main road (the remainder of that duration is the driver's {\em safety margin} to cross the road). Naturally, this further improves the capacity of the minor road, but our main conclusion that capacity is determined by the {\em entire} headway distribution (and not by its mean only) remains equally valid. Another obvious extension concerns networks consisting of multiple intersections.

\section*{Acknowledgments} The research of Abhishek and M. Mandjes is partly funded by NWO Gravitation project {\sc Networks}, grant number 024.002.003. The authors thank Onno Boxma (Eindhoven University of Technology) and Bart van Arem (Delft University of Technology) for helpful discussions. 

\bibliographystyle{abbrv}

\end{document}